\def\fnote#1{\footnote}

\documentclass[12pt, amsfonts, amssymb, epsfig, verbatim]{article} 
\usepackage{epsfig}

\newcommand{\prf}{\noindent{\it {Proof.}}\ }
\def\squaresymb{\hbox{\vrule\vbox{\hrule\phantom{o}\hrule}\vrule}}
\newcommand{\prfend}{\nolinebreak[4] \nopagebreak[4] \hfill $\squaresymb$}

\newtheorem{pro}{Proposition}[section]
\newtheorem{thm}[pro]{Theorem}
\newtheorem{lem}[pro]{Lemma}
\newtheorem{cor}[pro]{Corollary}

\newtheorem{defi}[pro]{Definition}

\newcommand{\circledast}{\oslash}
\newcommand{\HH}{\mbox{\sf%
\makebox[0ex][l]{\rule[-.025ex]{.15ex}{1.55ex}}%
\makebox[0ex][l]{\rule[1.53ex]{.4ex}{.05ex}}%
\makebox[.05ex][l]{\rule[-.025ex]{.4ex}{.05ex}}H}}

\def\cstar{${}$\circledast${}$}
\newcommand{\td}{\downarrow}
\newcommand{\tD}{\Downarrow}

\newcommand{\la}{\langle}
\newcommand{\ra}{\rangle}

\newcommand{\cA}{{\cal A}}
\newcommand{\cB}{{\cal B}}

\newcommand{\cG}{{\cal G}}
\newcommand{\cH}{{\cal H}}

\newcommand{\cP}{{\cal P}}
\newcommand{\cQ}{{\cal Q}}
\newcommand{\cR}{{\cal R}}

\newcommand{\bS}{{\bf S}}

\newcommand{\mcigp}{I_{\cG}}

\newcommand{\mcip}{Ip\;}
\newcommand{\bC}{{\bf C}}

\begin{document}
\hyphenation{de-com-pos-a-bi-li-ty}
\baselineskip=14pt

\centerline{{\Large {\bf Additive induced-hereditary properties and unique factorization }}}

\bigskip
\centerline{Grzegorz Arkit}
\smallskip
\centerline{{\it Institute of Mathematics}}
\centerline{{\it University of Zielona G\'ora}}
\centerline{{\it Podg\'orna 50, 65-246 Zielona G\'ora, Poland}}
\centerline{\scriptsize{e-mail:{\it G.Arkit@im.uz.zgora.pl}}}
\medskip
\centerline{Alastair Farrugia\footnote{Studies are fully funded by the Canadian government through a Canadian Commonwealth Scholarship.}}
\smallskip
\centerline{{\it Department of Combinatorics and Optimization }}
\centerline{{\it University of Waterloo}}
\centerline{{\it N2L 3G1, Canada}}
\centerline{\scriptsize{e-mail:{\it afarrugia@math.uwaterloo.ca}}}
\medskip
\centerline{Peter Mih\'ok\footnote{Research supported in part by
Slovak VEGA Grant 2/1131/21.}}
\smallskip
\centerline{{\it Department of Applied Mathematics}}
\centerline{{\it Faculty of Economics, Technical Univerzity}}
\centerline{{\it B.N{\v e}mcovej 32, 040 01 Ko{\v s}ice, Slovak republic}}
\centerline{{\it and}}
\centerline{{\it Mathematical Institute}}
\centerline{{\it Slovak Academy of Sciences}}
\centerline{{\it Gre\v s\'akova 6, 040 01 Ko{\v s}ice, Slovak republic}}
\centerline{\scriptsize{e-mail:{\it Peter.Mihok@tuke.sk}}}
\medskip
\centerline{Gabriel Semani\v sin\footnote{Research supported in part by
Slovak VEGA Grant 1/0424/03.}}
\smallskip
\centerline{{\it Institute of Mathematics}}
\centerline{{\it Faculty of Science, P.J. \v{S}af\'arik University}}
\centerline{{\it Jesenn\'a 5, 041 54 Ko\v{s}ice, Slovak republic}}
\centerline{\scriptsize{e-mail:{\it semanisin@science.upjs.sk}}}
\medskip
\centerline{Roman Vasky}
\smallskip
\centerline{{\it V.S.L. Software}}
\centerline{{\it Tr. SNP 61 - Kosil, 040 11 Ko{\v s}ice, Slovak republic}}
\centerline{\scriptsize{e-mail:{\it vasky@vsl.sk}}}
\medskip
\centerline{\today}

\begin{abstract}
We show that additive induced-hereditary properties of coloured 
hypergraphs can be uniquely factorised into irreducible factors.
Our constructions and proofs are so general that 
they can be used for arbitrary concrete categories of
combinatorial objects; we provide some examples of such 
combinatorial objects.


\smallskip
\noindent{\bf Primary Mathematics Subject Classification:} 05C15 \newline
\noindent{\bf Secondary Mathematics Subject Classification:} 18A10, 20L05, 05C65\newline
\noindent{{\bf Short title (35 characters):} Additive hereditary properties}

\end{abstract}

\section{Introduction}
Many problems treated in graph theory concern graph properties.
Roughly speaking, a graph property is a subset of the set of all
graphs, such as the family of planar graphs,
perfect graphs, interval graphs, claw-free graphs or hamiltonian
graphs. Some of these properties
have important common features that allow us to study them from 
a more general point of view.

{From} a  combinatorial aspect there is usually no need
to distinguish isomorphic copies of graphs and we therefore 
restrict our attention to unlabeled graphs. More precisely,
we require that a graph property be closed under isomorphism.

Many properties have the important feature of being 
also closed under taking some substructures. 
Consider a partial ordering $\preceq$ defined on the
set of graphs. 
A property $\cP$ is {\em $\preceq$-hereditary\/} if, whenever $G$ belongs to $\cal P$ and $H\preceq G$, then $H$ belongs to $\cal P$ as well.
Greenwell et al. proved in~\cite{grhe73} that such properties
are exactly those that can be characterized by the set of
$\preceq$-minimal forbidden substructures,
assuming that there is no infinite descending 
$\preceq$-chain of structures.
For example, the properties of $k$-degenerate graphs,
claw-free graphs and planar graphs are, respectively,
subgraph-hereditary, induced-subgraph-hereditary, and minor-hereditary.

Another important feature of many properties is that they are
closed under disjoint union of graphs. Such properties are said to be
{\em additive}. We show that this feature plays a substantial
role in the study of the structure of $\preceq$-hereditary
properties.

The language of graph properties can be successfully used to 
generalize ordinary vertex colouring. In a {\em proper colouring} 
each colour class must be an independent set. 
In so-called {\em generalized colouring}, each 
colour class must have a prescribed graph property.
Given a list of two or more properties, the class of all graphs 
that can be coloured  according to that list is said to be a {\em reducible} property.
One can immediately ask whether different lists correspond to different properties.
This is the {\em unique factorisation problem}, which was 
solved affirmatively for additive hereditary and additive 
induced-hereditary graph properties in~\cite{mise00, mi00, fari02a}.

In this paper we extend these results to induced-hereditary
properties of directed coloured hypergraphs. Moreover we show that our
result can be generalised beyond graphs and hypergraphs to
other combinatorial objects such as oriented graphs, 
or partially ordered sets. 

In Section 2 we introduce 
basic concepts and definitions that are used throughout the rest of the
paper.  In Section 3 we prove some necessary preliminary
results. Section 4 is devoted to canonical
factorisations of induced-hereditary properties of hypergraphs.
The Unique Factorisation Theorem for hypergraphs is presented in
Section 5. In the sixth section we introduce 
systems of objects of a concrete category, give some examples, and 
prove the Unique Factorisation Theorem for such systems.

\section{Basic concepts and definitions} 
In general we use standard graph and hypergraph terminology that
can be found, say, in~\cite{be76,be89}. 
For terminology related to hereditary properties of graphs and hypergraphs
we follow~\cite{bobr97}. In the next three sections we restrict our
attention to finite hypergraphs, without loops (hyperedges of size 1)
or multiple hyperedges. For the sake of brevity, we 
sometimes drop the ``hyper'' prefix, 
using ``edge'' instead of ``hyperedge''.

We will take our edges to be coloured and directed, so each
edge is not a set but an ordered tuple $(v_1, \ldots, v_r; c)$,
where the $v_i$'s are the vertices of the edge, and $c$ is its colour.
Isomorphisms must preserve the colour and direction of each edge.
The direction and colour actually make no difference, 
and are never mentioned in the proofs, so the reader might find it easier
to think about hypergraphs without colours or directions. We also point out in
advance that, if our properties contain only $k$-uniform hypergraphs, all our
constructions will only give $k$-uniform hypergraphs. Similarly, we may
restrict ourselves to hypergraphs with edge-colours taken from a prescribed set.

A~{\em hypergraph property} is any non-empty isomorphism-closed subclass of
hypergraphs. If $H$ belongs to a property $\cP$, then we also say that $H$ has property $\cP$.
The subhypergraph of $H$ induced by 
$U \subseteq V(H)$ is $H[U]$, with  edge-set 
$E(H[U]) := \{e \in E(H) | e \subseteq U\}$. $H'$ is an 
induced-subhypergraph of $H$ if it is  isomorphic to $H[U]$ 
for some $U \subseteq H$, and we write $H' \leq H$.

A property $\cP$ is {\em induced-hereditary}
if $H\in \cP$ implies that $K\in \cP$, for all $K\le H$.
A property is {\em additive} if it is closed under taking
disjoint union of hypergraphs. More precisely, if
$H_1=(V_{H_1},E_{H_1})$, $H_2=(V_{H_2},E_{H_2})$ are hypergraphs
with $V(H_1) \cap V(H_2) = \emptyset$,
then their disjoint union is the hypergraph $K=(V_{H_1}\cup
V_{H_2},E_{H_1}\cup E_{H_2})$.
A hypergraph is connected if and only if it cannot be expressed 
as a disjoint union of two hypergraphs.

Following the arguments in~\cite{bomi91} and \cite{ja02} 
one can easily verify that the set of all induced-hereditary
properties of hypergraphs ordered by set inclusion forms 
a completely distributive algebraic lattice,
which we shall denote by $\HH^a_{\le}$.
For many more details, applications and open problems concerning
hereditary and induced-hereditary properties we refer the reader
to~\cite{bobr97}.

Let $\cP_1,\cP_2, \dots, \cP_n$ be properties of hypergraphs.
A~$(\cP_1,\cP_2,\dots,\cP_n)$-{\em par\-ti\-tion} of a
hypergraph $H$ is a partition $(V_1,V_2,\dots,V_n)$ of the vertex set
$V(H)$ such that the induced subhypergraph $H[V_i]$ has property $\cP_i$, for
$i=1,2,\dots,n$. Note that $V_i$ could be empty for any $i$; equivalently, one can assume the null graph $K_0 = (\emptyset, \emptyset)$ to be contained in every property. 
If a hypergraph $H$ has a $(\cP_1,\cP_2,\dots,\cP_n)$-partition,
then we say that $H$ has property $\cP_1\circ\cP_2\circ\cdots\circ\cP_n$.
If $\cP_1=\cP_2=\cdots=\cP_n$ we simply write $\cP^n$ instead of
$\cP_1\circ\cP_2\circ\cdots\circ\cP_n$. 

Let $\cP$ be additive induced-hereditary; $\cP$ is {\em reducible} 
if there are additive induced-hereditary properties $\cP_1$ and $\cP_2$ 
such that $\cP=\cP_1\circ\cP_2$; otherwise, it is {\em irreducible}.
One may consider an alternative definition of reducibility in which $\cP_1$ and $\cP_2$ can be any two properties, not necessarily additive induced-hereditary. The two definitions turn out to be equivalent, but the proof of this non-trivial fact depends on the Unique Factorisation Theorem, and a further result characterising the existence of uniquely colourable graphs~\cite{brbu99, brbu02}, so we will stick with the first definition.

Unless stated otherwise, the properties we consider are additive induced-hereditary hypergraph properties. We will consider more general properties in the last section.

\section{Uniquely decomposable hypergraphs}
The main result of this section is the existence of uniquely $\cP$-decomposable hypergraphs, for every additive induced-hereditary property $\cP$. In fact, every hypergraph in $\cP$ is an induced-subhypergraph of a uniquely $\cP$-decomposable hypergraph.
\newline

Let $\cG$ be a set of hypergraphs. 
The  induced-hereditary property {\em generated by $\cG$\/} 
is $\la \cG\ra$, the smallest induced-hereditary property containing $\cG$. 
$\cG$ is a {\em generating set\/} for $\cP$ if $\la \cG\ra = \cP$. 
It is easy to see that:

\[\begin{array}{lll}
\la \cG\ra &=& \{G \mid \exists\, H \in \cG,\ G \leq H \}.
\end{array}\]

\noindent The {\em $*$-join\/} of $n$ hypergraphs $G_1, \ldots, G_n$ with disjoint vertex-sets is the set of all hypergraphs obtained by adding edges between the $G_i$'s; no new edges $e \subseteq V(G_i)$ are added:
\[ G_1 * \cdots * G_n := \{H \mid V(H) = \bigcup_{i=1}^n V(G_i), H[V(G_i)] = G_i\}. \] 
Given $n$ sets of hypergraphs, we define their $*$-join by 
\[ S_1 * \cdots * S_n := \bigcup \left(G_1 * \cdots * G_n\right)\] 
the union being over all ways of the selecting the $G_i$'s so that 
$G_i \in S_i$ for all $i$.
We note that this is just the same as $S_1 \circ \cdots \circ S_n$,
but it is aesthetically pleasing to have the $*$ notation. 

If $\cP_1, \ldots, \cP_n$ are additive properties, and $G_i \in \cP_i$ for all
$i$, then for all positive integers $k$ we have 
\[ kG_1 * \cdots * kG_n \subseteq \cP_1 \circ \cdots \circ \cP_n \] 
where $kG$ is the disjoint union of $k$ copies of $G$.
A {\em $\cP$-decomposition of $G$ with $n$ parts\/} is a partition
$(V_1,\ldots,V_n)$ of $V(G)$ such that for all $i$, $V_i \ne
\emptyset$, and for all positive integers $k$ we have $kG[V_1] *
\cdots * kG[V_n] \subseteq \cP$. The {\em $\cP$-decomposability
number\/} $dec_{\cP}(G)$ of $G$ is the maximum number of parts in a
$\cP$-decomposition of $G$; for $G \not\in \cP$ we put $dec_{\cP}(G) =0$. 
Thus $G$ is in $\cP$ if and only if $dec_{\cP}(G)\ge 1$.  
Also,  $G$ is {\em $\cP$-decomposable\/} if $dec_{\cP}(G) >1$. 
If $\cP$ is the product of two additive induced-hereditary properties, then
\emph{every} hypergraph in $\cP$ with at least two vertices is
$\cP$-decomposable. 

\begin{lem}
\label{ind-her-0}
{Let $\cP = \cP_1 \circ \cdots \circ \cP_m$, where the $\cP_i$'s are
additive properties. Then any $(\cP_1, \ldots, \cP_m)$-partition of a 
hypergraph $G$ is a $\cP$-decomposition of $G$. If
the $\cP_i$'s are induced-hereditary, then every hypergraph in $\cP$ with at least
$m$ vertices has a partition with all $m$ parts non-empty. 
}
\prfend
\end{lem}

A hypergraph $G$ is {\em $\cP$-strict\/} if $G$ is in $\cP$ but 
$G * K_1 \not\subseteq \cP$; we denote the set of $\cP$-strict hypergraphs by
$\bS(\cP)$.  
If $f(\cP) = \min\{|V(F)| \mid F \not\in \cP \}$, then $G * K_1 * \cdots
* K_1 \not\subseteq \cP$, where the $*$ operation is repeated $f(\cP)$
times. Thus, every $G \in \cP$ is an induced-subhypergraph of some
$\cP$-strict hypergraph (with fewer than $|V(G)| + f(\cP)$ vertices), and so
$\la \bS(\cP)\ra = \cP$. 
Similarly, $dec_{\cP}(G) < f(\cP)$.  

The {\em $\cP$-decomposability number\/} $dec_{\cP}(\cG)$ 
of a generating set $\cG$ of $\cP$
is $$\min\{dec_{\cP}(G) \mid  G \in \cG \};$$ the
{\em decomposability number\/} $dec(\cP)$ of $\cP$ is 
$dec_{\cP}(\bS(\cP))$. A property with $dec(\cP) = 1$ is
{\em indecomposable\/}. An indecomposable property is also irreducible
and it will turn out that the converse is also true. 

\begin{lem}
{Let $\cP_1, \ldots, \cP_m$ be induced-hereditary properties,
and let $G$ be a  $\cP_1 \circ \cdots \circ \cP_m$-strict
hypergraph. Then, for every $(\cP_1, \ldots, \cP_m)$-partition $(V_1,
\ldots, V_m)$ of $V(G)$, $G[V_i]$ is $\cP_i$-strict (and in particular
non-empty). 
}
\end{lem}

\prf
If $G[V_1]*K_1 \subseteq \cP_1$, then $G * K_1 \subseteq 
(G[V_1]*K_1) * G[V_2] * \cdots * G[V_m] \subseteq 
\cP_1 \circ \cdots \circ \cP_m$. 
\prfend
\newline

It follows that $dec(\cA \circ \cB) \geq dec(\cA) + dec(\cB)$, and
thus any factorisation of an additive induced-hereditary property $\cP$ has 
at most $dec(\cP)$ irreducible additive induced-hereditary factors.

\begin{lem}\ \cite{mise00}
\label{ind-her-2}
{Let $\cP$ be an induced-hereditary property and $G$ be a $\cP$-strict induced subhypergraph of $G' \in \cP$. Then $G'$ is $\cP$-strict, and
$dec_{\cP}(G) \geq dec_{\cP}(G')$. 
}
\end{lem}

\prf 
Every hypergraph in $G * K_1$ is an induced subhypergraph of a 
hypergraph in $G' * K_1$, 
so $G'$ must be $\cP$-strict. If $(V_1, \ldots, V_n)$ is a
$\cP$-decomposition of $G'$ with $n$ parts, then $(V_1 \cap V(G),
\ldots, V_n \cap V(G))$ is a $\cP$-decomposition of $G$; moreover, it
has $n$ parts unless, for some $i$, $V_i \cap V(G) = \emptyset$, which
is impossible because $G$ is $\cP$-strict. 
\prfend

\begin{lem}\ \cite{mise00}
\label{ind-her-3}
{If $\cG$ generates the induced-hereditary property $\cP$, 
then $dec_{\cP}(\cG) \leq dec_{\cP}(\bS(\cP))$, with equality if $\cG \subseteq
\bS(\cP)$.  
}
\prfend
\end{lem}

For $\cG \subseteq \cP$, and $H \in \cP$, let $\cG[H] := \{G \in \cG \mid
H \leq G \}$.  

\begin{lem}\ \cite{mise00}
\label{ind-her-4}
{Let $\cG$ generate the additive induced-hereditary property $\cP$,
and let $H$ be an arbitrary hypergraph in $\cP$. Then $\cG[H]$ also
generates $\cP$. 
}
\prfend
\end{lem}

For a generating set $\cG$, let $\cG^{\td} := \{G \in \cG \mid G \in
\bS(\cP),\ dec_{\cP}(G) = dec(\cP) \}$. The following is a simple
consequence of Lemmas \ref{ind-her-2} and \ref{ind-her-4}.

\begin{lem}\ \cite{mise00}
\label{ind-her-5}
{If $\cG$ generates the additive induced-hereditary property $\cP$,
then so does $\cG^{\td}$. 
}
\prfend
\end{lem}

A hypergraph $G$ is {\em uniquely $\cP$-decomposable\/} if it has 
exactly one $\cP$-dec\-omp\-osition with $dec_{\cP}(G)$ parts. 
Equivalently, $G$ is either $\cP$-indecomposable, or has exactly one $\cP$-decomposition with $n$ parts, for some $n \geq 2$; in the second case, $n$ must be $dec_{\cP}(G)$, as any decomposition with $n+1$ parts would give rise to ${n+1 \choose 2}$ decompositions with $n$ parts.

If $\cP =
\cP_1 \circ \cdots \circ \cP_n$, then by Lemma~\ref{ind-her-0} a
uniquely $\cP$-decomposable hypergraph $G$ with $dec_{\cP}(G) = n$ must be
uniquely $\{\cP_1, \ldots, \cP_n\}$-partition\-able (every $\{\cP_1,
\ldots, \cP_n\}$-partition gives the same unordered partition of
$V(G)$). 
If $(V_1, \ldots, V_n)$ is the unique $\cP$-decomposition of $G$, we
call the hypergraphs $G[V_1], \ldots, G[V_n]$ its {\em ind-parts} (although they
are themselves usually $\cP$-dec\-omposable). 

\begin{lem}
\label{ind-her-6}
{Let $\cP$ be an induced-hereditary property and let $G$ be a hypergraph in $\bS(\cP)$ with $dec_{\cP}(G) = dec(\cP)$, and
suppose that $G$ has a unique $\cP$-decomposition $(V_1, \ldots,
V_{dec(\cP)})$ with $dec(\cP)$ parts. If $G \leq H$, then $H \in
\bS(\cP)$, $dec_{\cP}(H) = dec(\cP)$, and, for any $\cP$-decomposition
$(W_1, \ldots, W_{dec(\cP)})$ of $H$, we can relabel the $W_i$'s so
that, for all $i$, $W_i \cap V(G) = V_i$. 
}
\prfend
\end{lem}

Let $d_0 = (U_1, U_2, \ldots, U_m)$
be a $\cP$-decomposition of a hypergraph $G$.
A $\cP$-decomposition $d_1=(V_1,V_2,\dots,V_n)$ of $G$ {\em respects\/} $d_0$ 
if no $V_i$ intersects two or more
$U_j$'s; that is, each $V_i$ is contained in some $U_j$, and so each
$U_j$ is a union of $V_i$'s. 

If $G$ is a hypergraph, then ${s}\cstar G$ denotes the set $G*G*\dots*G$, where
there are $s$ copies of $G$.  For $G^* \in {s}\cstar G$,
denote the copies 
of $G$ by $G^1, \ldots, G^s$. Then {\em $G^*$ respects $d_0$\/} if
$G^* \in s G[U_1] * \cdots * s G[U_m]$; that is, an edge that intersects
different $G^i$'s must also intersect different $U_j$'s. 
A $\cP$-decomposition 
$d = (V_1, \ldots, V_n)$ of $G^*$ {\em respects $d_0$
uniformly\/} if, for each
$V_i$, there is a $U_j$ such that, for every $G^k$, 
$V_i \cap V(G^k)\subseteq U_j$.  
The decomposition of $G^k$ induced by $d$
is denoted $d|G^k$.

If $G$ is uniquely $\cP$-decomposable, its ind-parts {\em respect
$d_0$\/} if its unique $\cP$-decomposition with $dec_{\cP}(G)$ parts
respects $d_0$. If $G^*$ is uniquely $\cP$-decompos\-able, its
ind-parts {\em respect $d_0$ uniformly\/} if: (a) for some $s$,
 $G^* \in {s}\cstar
G$; (b) $G^*$ respects $d_0$; and (c) $G^*$'s unique
$\cP$-decomposition with $dec_{\cP}(G^*)$ parts respects $d_0$
uniformly. 

\bigskip

\begin{figure}[htb]
\begin{center}
\input{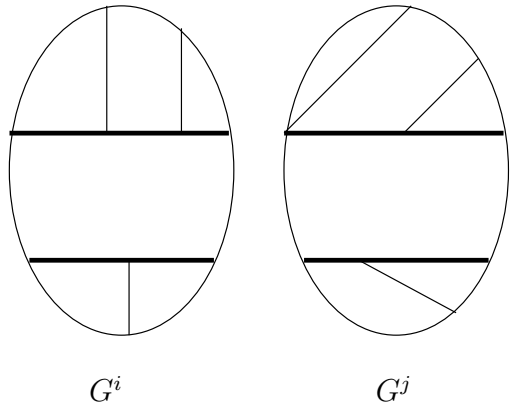} 
\caption{$d$ (vertical lines) respects $d_0$ (horizontal lines) uniformly}
\label{Fig-uniform-respect}
\end{center}
\end{figure}

\begin{figure}[htb]
\begin{center}
\input{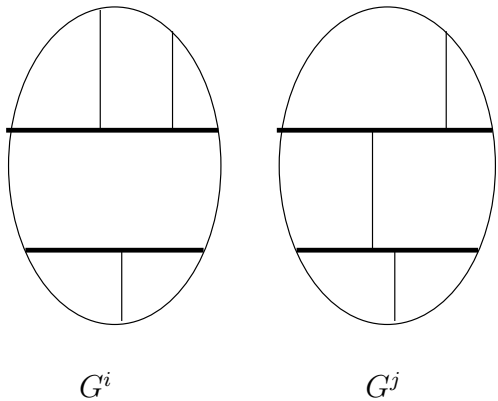} 
\caption{$d$ (vertical lines) respects $d_0$ (horizontal lines) on
both $G^i$ and $G^j$, but not uniformly} 
\label{Fig-respect}
\end{center}
\end{figure}

\bigskip

The {\em extension of $d_0$\/} to $G^*$ is the decomposition obtained
by repeating $d_0$ on each copy of $G$. If $G^*$ respects $d_0$, or if
it has a $\cP$-decomposition that respects $d_0$ uniformly, then the
extension of $d_0$ is also a $\cP$-decomposition of $G^*$. 
In particular, $G^*$ is a hypergraph in $\cP$.

We will sometimes write $G^i \cap U_x$ (or just $U_x$ when it is clear
we are referring to $G^i$) to mean the vertices of $G^i$ that
correspond to $U_x$, and $G^* \cap U_x$ (or just $U_x$, when it is
clear from the context) to mean $G^*[\bigcup_i (G^i \cap U_x)]$. 

The required result is a corollary of the following theorem of
Mih\'{o}k (see \cite{mi00}); he actually proved it when $m = n$
(Corollary~\ref{unique-super}), but very little modification is needed
to establish the general case, and we follow his proof and notation
rather closely. 
 
\begin{thm}
\label{m-resp-n}
{Let $G$ be a $\cP$-strict hypergraph with $dec_{\cP}(G) = n$, and let $d_0
= (U_1, U_2, \ldots, U_m)$ be a fixed $\cP$-decomposition of
$G$. Then there is a $\cP$-strict hypergraph $G^* \in {s}\cstar G$ (for some
$s$) that respects $d_0$, and moreover any $\cP$-decomposition of
$G^*$ with $n$ parts respects $d_0$ uniformly. 
}
\end{thm}

\prf
Let $d_i = (V_{i,1},V_{i,2}, \ldots ,V_{i,n}),\  i=1, \ldots, r$, be
the  $\cP$-decompositions of $G$ with $n$ parts which do not respect
$d_0$. Since $G$ is a finite hypergraph, $r$ is a nonnegative integer. If
$r=0$, take $G^* = G$; otherwise we will construct a  hypergraph $G^* =
G^*(r)\in s\cstar G$ as above, denoting the $s$ copies of $G$ by
$G^1, \ldots, G^s$. 

If the resulting $G^*$ has a $\cP$-decomposition $d$ with $n$ parts,
then, since $G$ is $\cP$-strict, $d|G^i$ will also have $n$ parts. The
aim of the construction is to add new edges $E^* = E^*(r)$ to
${s}G$ to exclude the possibility that $d|G^i = d_j$, for any $1 \leq
i \leq s, 1 \leq j \leq r$. 
Whenever we add an edge $e$, if $e$ intersects $G^i
\cap U_x$, it will also intersect some $G^j, i\not=j$, and some $U_y, x \not= y$; 
thus $G^*$ will respect $d_0$, and the hypergraphs 
constructed will always be in $\cP$.

We shall use two types of constructions.

\vspace{5mm}
{\noindent\bf Construction 1.} $G^i \Rightarrow G^j$.  

This is a
hypergraph in $2\cstar G$ such that, if $d$ is a $\cP$-decomposition of
$G^i\Rightarrow G^j$ and $d|G^i$ respects $d_0$, then $d|G^j$
respects $d_0$; moreover, $d$ respects $d_0$ uniformly on
$G^i\Rightarrow G^j$.

Since $G$ is $\cP$-strict, 
there is a hypergraph  $F \in (G * K_1)\setminus \cP$. 
Let $E'$ be the 
edges of $F$ that contain $z \in V(K_1)$. 
For $x=1,2,\dots,m$, 
let $E'_x$ be the set of edges from $E'$ that contain only $z$ and 
vertices of $U_x$, while $E'_{\overline{x}}$ is the set of edges 
from $E'$ that contain some vertex of $V(G) \setminus U_x$.
Let $G^i,G^j, i \not = j$, be disjoint copies of $G$; for every $x$, and 
every vertex $v \in U_x \cap V(G_j)$, we add the edges of $E'_{\overline{x}}$ 
(with $v$ taking the place of $z$, and $G_i$ taking the place of $G$).
Note that $G^i\Rightarrow G^j\in 2G[U_1]*2G[U_2]*\cdots
*2G[U_m]$.  Since $d_0$ is a $\cP$-decomposition of $G$, this implies
that $(G^i\Rightarrow G^j)\in \cP$.

Let $d=(V_1,V_2, \dots, V_\ell)$ be a $\cP$-decomposition of
$H=(G^i\Rightarrow G^j)$ such that $d|G^i$
respects $d_0$, but $d|G^j$ does not respect $d_0$ 
(or at least, not in the same manner, i.e., 
$d$ does not respect $d_0$ uniformly).
Then there exist $k$ and $x \ne y$ such that $V_k \cap G^i \subseteq U_y$, 
but some $v \in V_k \cap G^j$ belongs to $U_x$. 
We can add the edges corresponding to $E' \setminus E'_x$ 
because they contain at least one vertex $w$ of $G^i \cap U_x$
(so $w \not\in V_k$). But then, $F$ is an induced 
subhypergraph of a hypergraph in 
$ H[V_1]*H[V_2]*\dots*H[V_\ell]$, which implies $F\in\cP$, a contradiction.

\vspace{5mm}
{\noindent\bf Construction 2.} $m\bullet k_t G$.  

For a $\cP$-decomposition 
$d_t=(V_{t, 1},V_{t, 2},\dots,V_{t, dec_{\cP}(G)})$
 of $G$ that does not respect $d_0$,
$m\bullet k_t G$ is a
hypergraph in $(mk_t)\cstar G$ having no $\cP$-decomposition
$d=(W_1,W_2,\dots,W_{dec_{\cP}(G)})$ such that, for all of the $mk_t$
induced copies $G^i$ of $G$, $d|G^i = d_t$.

\begin{figure}[htb]
\begin{center}
\input{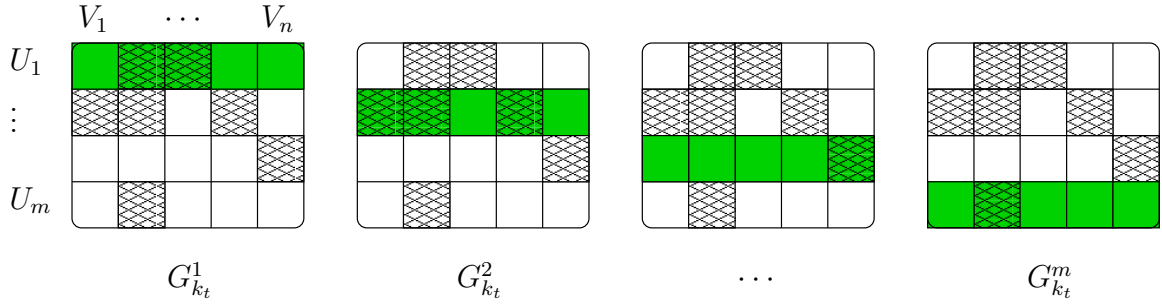} 
\caption{$m \bullet k_t G$ - we only put edges between the $m$ shaded parts}
\label{Fig-constr-2}
\end{center}
\end{figure}

Let $n=dec_{\cP}(G)$ and
let  $A_{i,j}(t)$ denote $U_i \cap V_{t,j}, 1 \leq i
\leq m, 1 \leq j \leq n$. Since $d_t$ does not respect $d_0$, at
least $n+1$ sets $A_{i,j}(t)$ are nonempty. Because $dec_{\cP}(G)=n$,
there exists a positive integer $k_t$ such that $k_t
G[A_{1,1}(t)] * k_t G[A_{1,2}(t)] * \ldots * k_t
G[A_{m,n}(t)] \not \subset \cP$. Fix a hypergraph  $F_t \in (k_t
G[A_{1,1}(t)] * k_t G[A_{1,2}(t)] * \ldots * k_t
G[A_{m,n}(t)])\setminus \cP$. 
Note that $F_t$ differs from $k_t G$ only in the edges that intersect at least two different $U_i$'s, or at least two different $V_j$'s.

The $U_i$'s form a $\cP$-decomposition of $k_t G$, so we can replace the edges of $k_t G$ that intersect at least two $U_i$'s, with the edges of $F_t$ that intersect at least two $U_i$'s, and still remain in $\cP$. If, in the resulting hypergraph $\tilde{H}$, the $V_j$'s also formed a $\cP$-decomposition, we could replace the edges of $\tilde{H}$ that intersect at least two different $V_j$'s with the edges of $F_t$ that intersect at least two different $V_j$'s, and still remain in $\cP$. But this is impossible because we would then have $F_t \in \cP$.

The only problem with $\tilde{H}$ is that, in order to construct it, we altered edges \emph{inside} the $k_t$ copies that we had of $G$.
We therefore construct $m \bullet k_t G$ by taking $m$ disjoint
copies of $H = k_t G$, denoted by $H^j, j=1,2, \ldots m$, and
adding edges between $H^1 \cap U_1, H^2 \cap U_2, \ldots, H^m \cap
U_m$. Specifically, suppose an edge of $F_t$ intersects $U_{a_1}, \ldots, U_{a_r}$ ($1 \leq a_1 < \cdots < a_r \leq m$, $r \geq 2$); then in $m \bullet k_t G$ we put a corresponding edge that intersects $H^{a_1} \cap U_{a_1}, \ldots, H^{a_r} \cap U_{a_r}$.

Suppose $d=(W_1,W_2,\dots,W_{dec_{\cP}(G)})$ is a $\cP$-decomposition
of $m\bullet k_t G$ such that, for every one of the $mk_t$
induced copies $G^i$ of $G$, $d|G^i=d_t$. Then $H^1 \cap U_1, \ldots, H^m \cap U_m$ induce a copy of the hypergraph $\tilde{H}$ from which we could obtain $F_t$ while still remaining in $\cP$, thus getting a contradiction as above.

\begin{figure}[htb]
\begin{center}
\input{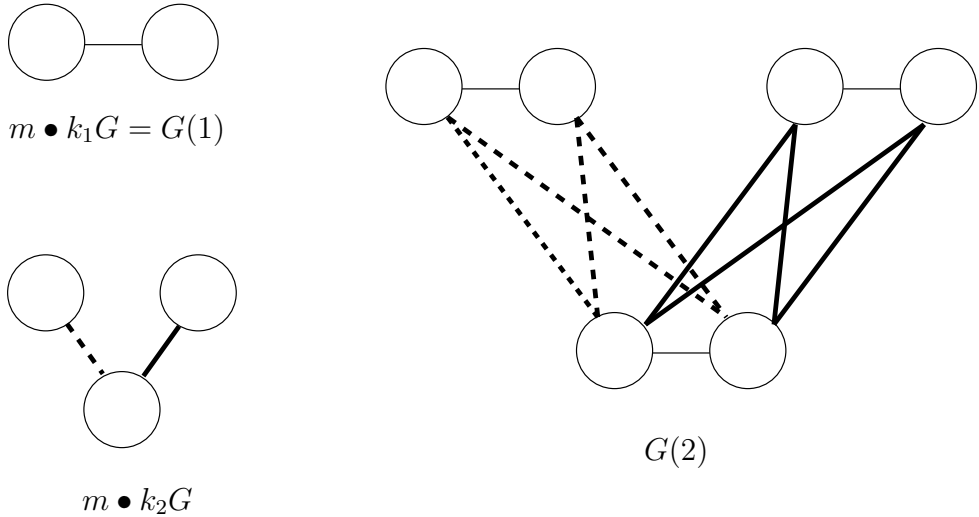} 
\caption{Constructing $G(2)$ from $G(1)$ and $m \bullet k_2 G$}
\label{Fig-constr-G-2}
\end{center}
\end{figure}

\vspace{5mm}
We now construct $G^*$ as follows.
First let $G(0) := G$ and $G(1) := m \bullet k_1 G$. For $1 < \ell \leq r$, construct
$G(\ell)$ by taking  $m k_\ell$ disjoint copies
$G(\ell-1)^1,\dots,G(\ell-1)^{mk_\ell}$ of $G(\ell-1)$.  For each copy 
of $G$ in $G(\ell-1)^i$ and each copy 
of $G$ in $G(\ell-1)^j$, we add the
edges between them that are between the $i^{{th}}$ and
$j^{{th}}$ copies of $G$ in $m\bullet k_\ell G$.  (See 
Figure~\ref{Fig-constr-G-2}.)

Finally, from $G(r)$, which is in, say, $s\cstar  G$, consisting of copies
$G^1$, $G^2$, \dots, $G^s$ of $G$, we create $G^*$ by adding 
two more copies $G^{+}$ and
$G^{-}$ of $G$.  
We add edges between $G^{+}$ and $G^{-}$ to create the hypergraph $G^{-}\Rightarrow G^{+}$, and, for each $i = 1,\dots,s$, we add edges 
to obtain $G^i\Rightarrow G^{-}$ and $G^{+}\Rightarrow G^i$.

Let $d$ be a $\cP$-decomposition of $G^*$ with $n$ parts (it might be
that none exists, in which case we are done).  
For $1 \leq \ell \leq r$, if every copy of $G(\ell-1)$ in $G(\ell)$ contains
a copy of $G$ for which $d|G = d_\ell$, then we would have $m k_\ell$ such
copies of $G$ inducing a copy of $m \bullet k_\ell G$, which we know is
impossible. So by induction from $r$ to $1$, there is a copy $G^p$
of $G$ for which $d|G^p$ is none of $d_1,d_2,\dots,d_r$.  Thus,
$d|G^p$ respects $d_0$.  But $G^p\Rightarrow G^{-}$ is an induced
subgraph of $G^*$, so $d|G^{-}=d_0$ (and in fact $d$ respects $d_0$
uniformly on these two copies of $G$).  Similarly, $d|G^{+}$ respects
$d_0$ and, in the same way, $d$ respects $d_0$ uniformly, as required.  
\prfend

\begin{cor}
\label{respect}
{Let $G$ be a $\cP$-strict hypergraph with $dec_{\cP}(G) = dec(\cP)$, and
let $d_0 = (U_1, U_2, \ldots,$ $U_m)$ be a fixed $\cP$-decomposition
of $G$. Then there is a $\cP$-decomposition of $G$ with exactly
$dec(\cP)$ parts that respects $d_0$. 
}
\end{cor}

\prf
In Theorem~\ref{m-resp-n}, since $G^* \geq G$ we know $G^*$ is
$\cP$-strict, and so $dec(\cP) \leq dec_{\cP}(G^*) \leq dec_{\cP}(G) =
dec(\cP)$. Thus $G^*$ has at least one $\cP$-decomposition $d$ with
$dec(\cP)$ parts; $d|G$ also has $dec(\cP)$ parts (since $G$ is
$\cP$-strict) and respects $d_0$. 
\prfend

\begin{cor}\ \cite{mise00}
\label{unique-super}
{Let $G$ be a $\cP$-strict hypergraph with $dec_{\cP}(G) = n$, and let $d_0
= (U_1, U_2, \ldots,$ $U_n)$ be a fixed $\cP$-decomposition of $G$ with
$n$ parts. Then there is a $\cP$-strict hypergraph $G^* \in{s}\cstar G$ (for
some $s$) which has a unique $\cP$-decomposition $d$ with $n$ parts,
and $d|G^j = d_0$ for all $j$. 
}
\end{cor}

\prf
The only $\cP$-decomposition of $G$ with $n$ parts that respects $d_0$
is $d_0$ itself (since here $d_0$ has exactly $n$ parts). Thus in
Theorem~\ref{m-resp-n}, the only possible decomposition of $G^*$ with
$n$ parts is the extension of $d_0$, which is a $\cP$-decomposition of
$G^*$ by construction. 
\prfend

\medskip

The set of $\cP$-strict,
uniquely $\cP$-decomposable hypergraphs with $dec_{\cP}(G) = dec(\cP)$ is
denoted $\bS^{\tD}(\cP)$, or just $\bS^{\tD}$. By
Lemma~\ref{ind-her-5} and  Corollary~\ref{unique-super} 
$\bS^{\tD}$ is a generating set for $\cP$; in fact, for any 
$G \in \bS^{\tD}$
and any specific $\cP$-decomposition $d$ of $G$, we can
find a hypergraph in $\bS^{\tD}$ that contains $G$ and whose ind-parts uniformly
respect $d$. 

\begin{cor}
\label{unique-respect}
{Let $G$ be a $\cP$-strict hypergraph with $dec_{\cP}(G) = dec(\cP)$, and
let $d_0 = (U_1, U_2, \ldots,$ $U_m)$ be a fixed $\cP$-decomposition
of $G$. Then there is a uniquely $\cP$-decomposable $\cP$-strict hypergraph
$G^* \geq G$ whose ind-parts respect $d_0$ uniformly. 
}
\prfend
\end{cor}

\section{Canonical factorisations}
A property $\cP$ is {\em indecomposable} if $dec(\cP) = 1$. In this section we show that every additive induced-hereditary property $\cP$ has a factorisation into $dec(\cP)$ additive induced-hereditary properties; this establishes the important fact that $\cP$ is irreducible iff it is indecomposable. In the next section we will show that whenever $\cP$ has a factorisation into indecomposable factors, there must be exactly $dec(\cP)$ of them, and this factorisation must be unique.

\begin{lem}\ \cite{sc86}
\label{ordered}
{A generating set $\cG' = \{G_1, G_2, \ldots\}$ for $\cP$ contains an ordered generating
set 
$\cG^\odot = \{G_{k_1}, G_{k_2}, \ldots\}\subseteq \cG'$ satisfying
$G_{k_1} \leq G_{k_2} \leq \cdots$\,.
}
\end{lem}

\prf
Since $\cP$ contains only finite graphs, it is countable, say $\cP = \{H_1, H_2, \ldots \}$. Pick $G_{k_1}$ arbitrarily. For each $i$, by additivity, $G_{k_i} \cup H_i$ is in $\cP$, so there is a $k_{i+1}$ for which $(G_{k_i} \cup H_i) \leq G_{k_{i+1}}$. \prfend

Recall that $\bS^{\tD} := \bS^{\tD}(\cP)$ is the set of uniquely $\cP$-decomposable hypergraphs with decomposability $n := dec(\cP)$. By Lemma~\ref{ind-her-5} and  Corollary~\ref{unique-super}, $\bS^{\tD}$ is a generating set for $\cP$. 
By Lemma~\ref{ordered} there is an ordered generating set $\cG \subseteq \bS^{\tD}$ for $\cP$. 

For a hypergraph $G \in \bS^{\tD}$ whose unique $\cP$-decomposition is $(V_1, \ldots, V_n)$, the set of ind-parts is $\mcip(G) := \{G[V_1], \ldots, G[V_n]\}$. The set of all ind-parts from $\cG$ is $\mcigp := \bigcup (\mcip(G): G \in \cG)$.
For $F \in \mcigp$ and $G \in \cG$, $m(F,G)$ is the {\em multiplicity} of $F$ in $G$: the number of different (possibly isomorphic) ind-parts of $G$ which contain $F$ as an induced-subhypergraph. The multiplicity of $F$ in $\cG$ is $m(F)=\max\{m(F,G) \mid G\in\cG\}$; clearly $1\leq m(F)\leq n = dec({\cP})$.

\begin{lem} \label{tech_lem} 
Let $G, H,$ be hypergraphs in $\bS^{\tD}$. If $G \leq H$, then each ind-part of $G$ is an induced-subhypergraph of a distinct ind-part of $H$.
\end{lem} 

\prf
Let the ind-parts of $G$ and $H$ be $(G_1, G_2, \ldots G_n)$ and $(H_1, \ldots, H_n)$, respectively. $(H_1 \cap G, H_2 \cap G, \ldots, H_n \cap G)$ is a $\cP$-decomposition of $G$, where we use $H_k \cap G$ to denote $V(H_k) \cap V(G)$. If $H_k \cap G = \emptyset$ for some $k$, then $G$ would not be $\cP$-strict, a contradiction. So we have a $\cP$-decomposition of $G$ with $n$ parts. Because $G$ is in $\cG$, there is only one such decomposition, so without loss of generality $H_k \cap G = V(G_k), \ k = 1,\ldots, n$.
\prfend

For convenience, we will talk of the hypergraph induced by ind-parts $G_1, G_2, \ldots,$ when we actually mean the subhypergraph induced by $V(G_1) \cup V(G_2) \cup \cdots$\,.

\begin{thm}\ \cite{mise00}
\label{ind-her-factorn}
{An additive induced-hereditary property $\cP$ has a factorisation into $dec(\cP)$ (necessarily indecomposable) additive induced-hereditary factors.
}
\end{thm}

\prf
We proceed by induction on $dec(\cP)$.
If $dec(\cP) = 1$ there is nothing to do. So let every hypergraph $G \in \bS(\cP)$ with at least two vertices be $\cP$-decomposable. We will either factorise at once into $n := dec(\cP)$ properties, or into properties $\cQ, \cR$ such that $dec(\cP) = dec(\cQ) + dec(\cR)$.

\medskip

\noindent{\bf Case 1.} $m(F) = k < dec(\cP)$, for some $F \in \mcigp$. 

Let $G\in\cG$ be a generator of
$\cP$ for which $m(F,G)=k$. By Lemma~\ref{ind-her-4}, $\cG[G]$ generates $\cP$; by Lemma \ref{tech_lem}, for every
generator $H\in\cG[G]$, $m(F,H)=k$, so $\cG_{F} := \{G' \in \cG \mid m(F,G') = k\}$ is a generating set. For $H \in \cG_F$, let $H_F$ be the subgraph induced by the $k$ ind-parts which contain $F$, and $H_{\overline F}$ the subhypergraph induced by the $n-k$ other ind-parts.
Let the induced-hereditary properties
$\cQ_{F}$ and $\cQ_{\overline F}$ be generated by 
$\{H_F \mid H \in \cG_F\}$ and $\{H_{\overline{F}} \mid H \in \cG_F\}$, respectively.

We claim that $\cP =\cQ_{F}\circ\cQ_{\overline{F}}$. It is easy to see that
$\cP\subseteq\cQ_{F}\circ\cQ_{\overline{F}}$. Conversely, let $H$ be in $\cQ_F\circ\cQ_{\overline{F}}$. Then $H \in H^1_F \ast H^2_{\overline F}$, for some
$H^1, H^2 \in \cG_F$. Let $H'$ be a hypergraph in $\cG$ such that $H^1 \cup H^2 \leq H'$. By Lemma~\ref{tech_lem}, and because the maximum multiplicity of $F$ in $\cG$ is $k$, $H^1_F \leq H'_F$ and $H^2_{\overline F} \leq H'_{\overline F}$. Since $H'_F \ast H'_{\overline F} \subseteq \cP$, we have $H^1_F \ast H^2_{\overline F} \subseteq \cP$, implying 
$H\in\cP$. Hence $\cP =\cQ_F\circ\cQ_{\overline F}$.

To establish additivity of $\cQ_F$, consider
$G_F, H_F \in \cQ_F$, for some $G, H \in \cG_F$. 
Because $\cG_F$ generates $\cP$, there is some $L \in \cG_F$ such that
$(G \cup H) \leq L$. By Lemma~\ref{tech_lem}, $G_F \leq L_F$ and 
$H_F \leq L_F$, so $(G_F \cup H_F) \leq L_F \in \cQ_F$, and thus $(G_F \cup H_F) \in \cQ_F$. Additivity of $\cQ_{\overline F}$ is proved similarly.


Finally, we want to show that every $H_F \in \cG_F$ has $\cQ_F$-{de\-com\-pos\-a\-bi\-li\-ty} at least $k$, and every $H_{\overline F} \in \cG_{\overline F}$ has $\cQ_{\overline F}$-{de\-com\-pos\-a\-bi\-li\-ty} at least $n-k$. This will imply that $n = dec(\cP) \geq dec(\cQ_F) + dec(\cQ_{\overline F}) \geq k + (n-k) = n$, and thus $dec(\cQ_F) = k$ and $dec(\cQ_{\overline F}) = n-k$. Since $k < n$ and $n-k < n$, the factorisation result will follow by induction.

So consider $H_F \in \cG_F$, and let $H'_F$ be in $H_1 * \cdots *  H_k$, where $H_1, \ldots, H_k$ are the ind-parts of $H$ that contain $F$. 
Consider $H'' := H \cup H'_F$; since $H \leq H''$, $H''$ is $\cP$-strict, and $dec_{\cP}(H'') = dec(\cP)$; moreover, $H''$ has a $\cP$-decomposition where the parts are $2H_1, \ldots, 2H_k, H_{k+1}, \ldots, H_n$ . By Corollary~\ref{unique-super}, there is a uniquely $\cP$-decomposable graph $H^* \in s \cstar H''$ whose 
ind-parts are $2sH_1, \ldots, 2sH_k, sH_{k+1}, \ldots, sH_n$. Thus $H'_F \leq H^*_F$, and so $H'_F \in \cQ_F$. Since $H'_F$ was arbitrary, $H_1, \ldots, H_k$ give a $\cQ_F$-decomposition of $H_F$, and so $dec_{\cQ_F}(H_F) \geq k$ as required.
The proof that $dec(\cQ_{\overline F}) \geq n-k$ is similar.
\medskip

\noindent{\bf Case 2.} $m(F)=n:=dec({\cP})\geq 2$ for each $F\in \mcigp$. 

Let $\cQ$ be the induced-hereditary property generated by $\mcigp$. It is
easy to see that $\cP \subseteq \cQ^n$.
The converse inclusion, $\cQ^n \subseteq \cP$, and the additivity and indecomposability of $\cQ$, follow as in Case 1.
\prfend

\begin{cor}\ \cite{mise00}
\label{ind-her-irr-ind}
{An additive induced-hereditary property is irreducible if and only if it is indecomposable.
}
\prfend
\end{cor}

\section{Unique factorization theorem for hypergraphs}
\label{sec-ind-proof}
To prove unique factorisation, we shall first show that 
the number of factors must be exactly $dec(\cP)$, and 
then show that any two factorisations with $dec(\cP)$ 
factors must be the same.
%
%
In the case of additive hereditary properties, there is a simple direct proof of the following result implicit in~\cite[Lemma 2.1]{mi00}.

\begin{lem}
\label{unique-respect-2}
{
Let $\cP$ be $\cQ_1 \circ \cdots \circ \cQ_m$. 
Let $G$ be a $\cP$-strict, uniquely
$\cP$-decomposable graph with $dec_{\cP}(G) = dec(\cP)$, 
and let $(W_1, \ldots, W_m)$ be a 
$(\cQ_1, \ldots, \cQ_m)$-partition of $G$. 
Then each $W_j$ is a union of ind-parts of $G$.
}
\end{lem}

\prf
By Corollary~\ref{unique-respect} there is a a uniquely
$\cP$-decomposable graph $G^* \geq G$ whose ind-parts respect
the $W_j$'s. Now the ind-parts of $G$ are just the restriction
of the ind-parts of $G^*$.
\prfend

\begin{thm}
\label{ind-her-m=n-ver-2}
{Let 
$\cQ_1 \circ \cdots \circ \cQ_m$ be a factorisation of the
additive induced-hereditary property $\cP$ into indecomposable additive
induced-hereditary properties. Then $m = dec(\cP)$. 
}
\end{thm}
 
\prf
By Lemma~\ref{ind-her-0} any $\cP$-strict graph $G$ has $dec_{\cP}(G) \geq m$, so $dec(\cP) \geq m$. 
To prove the reverse inequality, note that $\cP$ is generated by the set of $\cP$-strict, uniquely $\cP$-decomposable graphs with minimum decomposability, by Lemma~\ref{ind-her-5} and Corollary~\ref{unique-super}. 
This contains an ordered generating set $\cG$, say $G_1 \leq G_2 \leq \cdots,$ as constructed in Lemma~\ref{ordered}. So we have

\noindent (a) each $G_r$ is $\cP$-strict and uniquely $\cP$-decomposable, 
with $dec_{\cP}(G_r) = n$.

Let $\cP_1 \circ \cdots \circ \cP_n$ be a factorisation of $\cP$, relative to $\cG$, into $n:= dec(\cP)$ indecomposable factors, as constructed in Theorem~\ref{ind-her-factorn}. Then we can label the ind-parts of each $G_r$ as $G_{1,r}, \ldots, G_{n,r}$, so that:

\noindent (b) for each $i$, the $G_{i,r}$'s are ordered by inclusion, 
say $G_{i,1} \leq G_{i,2} \leq \cdots,$ and they form a generating set for $\cP_i$.

The first part of (b) follows from Lemma~\ref{tech_lem} and the fact that the $G_r$'s are themselves ordered. The second part follows from the proof of Theorem~\ref{ind-her-factorn};  in Case 2 it is clear. 
When $m(F) = k < n$ for some $F$, let $m(F, G_s) = k$, with $F$ contained in, say, $G_{1,s}, \ldots, G_{k,s}$; then $m(F, G_t) = k$ for all $t \geq s$, and $G_{1,t}, \ldots, G_{k,t}$ are the ind-parts in $\cQ_F$. We remove $G_1, \ldots, G_{s-1}$ from the generating set $\cG$, and assertion (b) then follows by induction on $dec(\cP)$.
\newline

For each $i$ and $j$ take an arbitrary $X_{i,j}\in \cP_i\setminus \cP_j$; 
if $\cP_i\setminus \cP_j=\emptyset$, then set $X_{i,j}$ to be the null graph $K_0$.
We set $\displaystyle H_{i}:=\bigcup_{j}X_{i,j}$; note that $H_{i}$ is in $\cP_i$.
The important point is that if $\{Y_1,Y_2,\dots,Y_n\}$ is an
unordered $(\cP_1,\dots,\cP_n)$-partition of some graph 
$G$ such that, for each $i=1,2,\dots,n$, 
$H_{i}\leq G[Y_i]$, then $G[Y_i]\in\cP_i$. If not, let $k_1, \ldots, k_r$
be the indices for which $G[Y_{k_j}] \not\in \cP_{k_j}$; then there is a permutation
$\varphi$ of the $k_j$'s such that $G[Y_{k_j}] \in \cP_{\varphi(k_j)} \not= \cP_{k_j}$, and we get a contradiction when we consider any $\cP_{k_s}$ that is inclusion-wise maximal among the $\cP_{k_j}$'s.
We must have $H_i \leq G_{i, r}$ for $r$ sufficiently large, so we can omit finitely many $G_r$'s to get:

\noindent (c) $H_i \leq G_{i,1}$ for each $i$.

Properties (a, b, c) guarantee that $(G_{1,r}, \ldots, G_{n,r})$ is the unique \emph{ordered}
$(\cP_1, \ldots, \cP_n)$-partition of $G_r$. For each $G_r$ we fix some ordered
$(\cQ_1, \ldots, \cQ_m)$-partition (it must have at least one such partition).
By Lemma~\ref{unique-respect-2} each $\cQ_i$-part is the union of ind-parts of $G_r$, that is,
there is a partition $(S_{1,r},\ldots,S_{m,r})$ 
of $\{1, \ldots, n\}$ such that $G_r[\cup_{s \in S_{j,r}} V(G_{s,r})]  \in \cQ_j$,
for each $j = 1, \ldots, m$.

By (b), the partition $(S_{1,r},\ldots,S_{m,r})$ also works for $G_1, G_2, \ldots, G_{r-1}$. Since there are only finitely many partitions of $\{1, \ldots, n\}$, one of them must
appear infinitely often, so we can use this partition of the ind-parts for all $r$; let
it be $(S_1, \ldots, S_m)$.

We want to prove that $\cQ_1 = \prod_{s \in S_1} \cP_s$. Since $\cQ_1$ is 
irreducible, this will imply that $|S_1| = 1$; the same reasoning applies
to $S_j, j = 2, \ldots, m$, so that we must have $m = n$.

Without loss of generality, $S_1 = \{1, \ldots, q\}$. 
Let $A$ be a graph in $\cQ_1$. Note that in $G_1$, the ind-parts
$G_{1,1}, \ldots, G_{q,1}$ form a graph in $\cQ_1$.
Let $v$ be a vertex of $G_{1,1}$, and let $N$ be the set of neighbours of $v$
in $G_{q+1,1}, \ldots, G_{n,1}$. Let $C$ be the graph formed from $G_1 \cup A$
by adding all possible edges between $A$ and $N$. $C$ has a
$(\cQ_1, \ldots, \cQ_m)$-partition with $V(A)$ in the $\cQ_1$-part, so it
is in $\cP = \cQ_1 \circ \cdots \circ \cQ_m$ and thus has a 
$(\cP_1, \ldots, \cP_n)$-partition.
Now $(G_{1,1}, \ldots, G_{n,1})$ is the unique $(\cP_1, \ldots, \cP_n)$-partition 
of $G_1$. If any vertex of $A$ is in $\cP_j$, $j > q$, then we could have put $v$
in $\cP_j$, a contradiction, so $A \in \cP_1 \circ \cdots \circ \cP_q$.

The reverse containment is proved similarly, but requires a bit more work.
Let $B$ be a graph in $\cP_1 \circ \cdots \circ \cP_q$, with 
$(\cP_1, \ldots, \cP_q)$-partition $(B_1, \ldots, B_q)$. 
We first create a graph $B' \in \cP$ that consists of several copies of $B$: for every tuple $(i_1, \ldots, i_q)$ such that $\cP_1 = \cP_{i_1}, \cP_2 = \cP_{i_2}, \ldots, \cP_q = \cP_{i_q}$, we put a copy of $B$ with $B_1, \ldots, B_q$ in the $\cP_{i_1}, \ldots, \cP_{i_q}$ part, respectively. This has an obvious $(\cP_1, \ldots, \cP_n)$-partition, say $(B'_1, \ldots, B'_n)$; note that $B'_j$ is empty iff $\cP_j$ is not equal to any of $\cP_1, \ldots, \cP_q$. 

As before, we take a vertex $v_1 \in G_{1,1}$ and let $N_{\overline 1}$ be $N(v_1) \setminus V(G_{1,1})$. Similarly, we take $v_2 \in G_{1,2}$ and let $N_{\overline 2}$ be $N(v_2) \setminus V(G_{1,2}$; and so on for $v_3, \ldots, v_n$ and 
$N_{\overline 3}, \ldots, N_{\overline n}$. Let $D$ be the graph formed from $G_1 \cup B'$ by adding all possible edges between $B'_1$ and $N_{\overline 1}$, $B'_2$ and $N_{\overline 2}$, \ldots. Then $D$ has a unique $(\cP_1, \ldots, \cP_n)$-partition, \emph{up to relabeling of identical properties}. 

Let $D$ be contained in $G_r$, for some $r$. By construction of $B'$, no matter which properties get labeled as $\cP_1, \ldots, \cP_q$, there will be a copy of $B$ contained in $G_r[V(G_{1,r}) \cup \cdots \cup V(G_{q,r})]$. This subgraph is in $\cQ_1$, so we are done.
\prfend
\newline

We will use the following construction of a generating set for $\cP$ to prove
unique factorisation. Suppose we are given a
factorisation $\cP = \cP_1 \circ \cdots \circ \cP_m$ into
indecomposable additive induced-hereditary factors, and, for each $i$,
we are given a generating set $\cG_i$ of $\cP_i$ and a graph $H_i \in
\cP_i$. By 
Lemmas~\ref{ind-her-4} and~\ref{ind-her-5}, the set $\cG_i^{\td}[H_i]
:= \{G \in (\cG_i \cap \bS(\cP_i)) \mid H_i \leq G,\  dec_{\cP_i}(G) = 1 \}$
is also a generating set for $\cP_i$. 
The $*$-join of these $m$ sets is then a generating set for $\cP$,
and even if we pick out just those graphs that are strict and
have minimum decomposability, we still have a generating set:

\begin{eqnarray*}
(\cG_1[H_1] * \cdots * \cG_m[H_m])^{\td}  :=  \{G' \in \bS(\cP)\ |\
dec_{\cP}(G') = dec(\cP), \textrm{ and } \forall\, i, \\ 
1\le i\le m,\ \exists\, G_i \in \cG_i^{\td}[H_i],\ G' \in G_1 *
\cdots * G_m\}. 
\end{eqnarray*} 

\begin{thm}
\label{ind-her-n}
{An additive induced-hereditary property $\cP$ can have only one
factorisation with exactly $dec(\cP)$ indecomposable factors. 
}
\end{thm}

\prf
Let $\cP_1 \circ \cdots \circ \cP_n = \cQ_1 \circ \cdots \circ \cQ_n$ 
be two factorisations of $\cP$ into $n := dec(\cP)$
indecomposable factors. 
Label the $\cP_i$'s inductively, beginning with $i=n$, so that, for each
$i$, $\cP_i$ is inclusion-wise maximal among
$\cP_1,\cP_2,\dots,\cP_i$.  For each $i,j$ such that $i>j$, if
$\cP_i\setminus \cP_j\ne\emptyset$, then let $X_{i,j}\in
\cP_i\setminus \cP_j$; if $\cP_i\setminus \cP_j=\emptyset$, then
$\cP_i= \cP_j$ and we set $X_{i,j}$ to be the null graph.  For each
$i$, set $\displaystyle H_{i,0}:=\bigcup_{j<i}X_{i,j}$.  Note
$H_{i,0}\in\cP_i$.  The important point is
 that if $\{Y_1,Y_2,\dots,Y_n\}$ is an
unordered $(\cP_1,\dots,\cP_n)$-partition of some graph 
$G$ such that, for each $i=1,2,\dots,n$, 
$H_{i,0}\leq G[Y_i]$, then, by reverse 
induction on $i$ starting at $n$, $G[Y_i]\in\cP_i$.

For each $i$, let $\cG_i = \{G_{i,0}, G_{i,1}, G_{i,2}, \ldots \}$ be a
generating set for $\cP_i$. We will construct another generating set
for each $\cP_i$ that will turn out to be contained in some $\cQ_j$;
for graphs $G_{i,s}, H_{i,s}$, we will use the second subscript to
denote which step of our construction we are in. 

For each $s \geq 0$, choose a graph $H'_{s+1} \in (\cG_1[H_{1,s},
G_{1,s}] * \cdots * \cG_n[H_{n,s}, G_{n,s}])^{\td}$, and find an
induced supergraph $H_{s+1}$ whose unique $\cP$-decomposition with
$dec(\cP)$ parts uniformly respects the obvious decomposition of
$H'_{s+1}$. We label as $H_{i,s+1}$ the ind-part of $H_{s+1}$ that
contains the graph from $\cG_i[H_{i,s}, G_{i,s}]$. Then, for each $i$,
$H_{i,0} \leq H_{i,1} \leq H_{i,2} \leq \cdots$ 

For $\cG_i[H_{i,s}, G_{i,s}]$ to be non-empty, we must have $H_{i,s}
\in \cP_i$. We know that the $H_{i,s+1}$'s give an unordered $\{\cP_1,
\ldots, \cP_n\}$-partition of $H_{s+1}$. From the earlier remark, for
$i=1,2,\dots,n$, 
$H_{i,s+1}\in \cP_i$.

The ind-parts of $H_s$ also form its unique $\{\cQ_1, \ldots,
\cQ_n\}$-partition.  Thus, there is some permutation $\varphi_s$ of
$\{1,2,\ldots, n\}$ such that, for each $i$, $H_{i,s} \in
\cQ_{\varphi_s(i)}$. Since there are only finitely many permutations
of $\{1,2,\ldots,n\}$, there must be some permutation $\varphi$ that
appears infinitely often. Now whenever $\varphi_t = \varphi$, we
have $H_{i,1} \leq H_{i,2} \leq \cdots \leq H_{i,t} \in
\cQ_{\varphi(i)}$ so by induced-heredity, for every $s\le t$,
$H_{i,s}$ is in $\cQ_{\varphi(i)}$.  Therefore,  we can take
$\varphi_s = \varphi$, for all $s$. By re-labelling the $\cQ_i$'s, we
can assume $\varphi$ is the identity permutation, so that $H_{i,s} \in
\cQ_i$ for all $i$ and $s$.  

Now for each $i$ and $s$, $G_{i,s-1} \leq H_{i,s}$, so that $\cH_i :=
\{H_{i,1}, H_{i,2}, \ldots \}$ is a generating set for $\cP_i$. But
$\cH_i \subseteq \cQ_i$, so $\cP_i = \la \cH_i\ra \subseteq \cQ_i$. 

By the same reasoning there is a permutation $\tau$ such that $\cQ_i
\subseteq \cP_{\tau(i)}$. We cannot relabel the $\cP_i$'s as well, but
if $\tau^k(i) = i$, then we have $\cP_i \subseteq \cQ_i \subseteq
\cP_{\tau(i)} \subseteq \cQ_{\tau(i)} \subseteq \cP_{\tau^2(i)}
\subseteq \cQ_{\tau^2(i)} \subseteq \cdots \subseteq \cP_{\tau^k(i)} =
\cP_i$, so we must have equality throughout; in particular, $\cP_i =
\cQ_i$ for each $i$. 
\prfend

\begin{thm}
\label{ind-her-uft}
{An additive induced-hereditary property $\cP$ has a unique factorisation
into irre\-ducible additive induced-here\-di\-tary factors, and the number of
factors is exactly $dec(\cP)$. }
{\nolinebreak[4] \nopagebreak[4] $\squaresymb$}
\end{thm}

\section{Unique factorization theorem for systems}

In this section we will present a common generalization of graphs, hypergraphs, digraphs and other combinatorial systems. We will use the basic elementary notions of category theory (see~\cite{pi91}) and 
deal only with concrete categories. 
A concrete category $\bC$  is a collection of {\it objects} and {\it arrows} called {\it morphisms}. An object in a concrete category $\bC$ is ``a set with structure''. We will denote the {\it ground-set} of the object $A$ by $V(A)$. The morphism between two objects is a ``structure preserving mapping''. Obviously, the morphisms of $\bC$ have to satisfy the axioms of the category theory (see e.g.~\cite{pi91}, page 1). The natural examples of concrete categories are: ${\bf Set}$ of sets, ${\bf FinSet}$ of finite sets, ${\bf Graph}$ of graphs, ${\bf Grp}$ of groups, ${\bf Poset}$ of partially ordered sets with structure preserving mappings, called homomorphisms of corresponding structures. In our investigations here we will need to consider 
{\it isomorphisms} i.e. structure preserving bijections between the ground-sets of objects only.

A simple finite hypergraph $H = (V,E)$ can be considered as a system of its hyperedges $E = \{e_1,e_2,\dots, e_m \}$, where edges are finite sets and the set of its vertices $V(H)$ 
is a superset of the union of hyperedges, i.e. $V \supseteq \bigcup_{i=1}^m e_i$. The following definition gives a natural generalization of hypergraphs or ``set-systems''.

\begin{defi}
Let $\bC$ be a concrete category. A {\em simple system of objects} of $\bC$ is an ordered pair
$S = (V, E)$, where $E = \{A_1,A_2,\dots, A_m \}$ is a finite set of the objects of $\bC$, such that the ground-set $V(A_i)$ of each object $A_i \in E$ is a finite set with at least two elements (i.e. there are no loops) and $V \supseteq \bigcup_{i=1}^m V(A_i)$.
\end{defi}

For example, graphs can be viewed as systems of objects of a concrete category of two-element sets with bijections as arrows, digraphs as special systems of objects of the category of posets, etc.

To generalize the proof of Unique factorization for coloured hypergraphs to arbitrary simple systems of objects 
(or shortly systems) 
we need to define ``isomorphism of systems'', ``disjoint union of systems'' and ``induced-subsystems'', respectively. We can do this in a natural way:

Let $S_1 = (V_1, E_1)$ and $S_2 = (V_2, E_2)$ be two simple systems of objects of a given concrete category $\bC$. 

The systems $S_1$ and $S_2$ are said to be isomorphic if 
there are two bijections:
$$ \phi : V_1 \longleftrightarrow V_2; \hspace{3cm} \psi : E_1 \longleftrightarrow E_2 ,$$
such that if $\psi(A_{1i}) = A_{2j}$ then $\phi / V(A_{1i}) : V(A_{1i}) \longleftrightarrow V(A_{2j})$ is an isomorphism of the objects $A_{1i} \in E_1$ and $A_{2j} \in E_2$ in the category $\bC$.

The disjoint union of the systems $S_1$ and $S_2$ is the system $S_1 \cup S_2 = ( V_1 \cup
V_2, E_1 \cup E_2)$, where we assume that $V_1 \cap V_2 = \emptyset$. 

A system is said to be connected if it cannot be expressed as a disjoint union of two systems.

The subsystem of $S_1$ induced by the set 
$U \subseteq V(S_1)$ is $S_1[U]$, with  objects 
$E(S_1[U]) := \{A_{1i} \in E(S_1) | V(A_{1i}) \subseteq U\}$. $S_2$ is an 
induced-subsystem of $S_1$ if it is  isomorphic to $S_1[U]$ 
for some $U \subseteq V(S_1)$. 

Using these definitions we can say, analogously as for hypergraphs, that an 
additive induced-hereditary property of simple systems of objects of a category $\bC$ is any class of systems closed under taking induced-subsystems, disjoint union of systems and isomorphism, respectively.
To prove the Unique Factorization Theorem for induced hereditary and additive properties of simple systems of objects of a concrete category $\bC$ we can follow the notions and constructions given in the previous Sections with some additional technical details, which we will omit here.

\ifx\undefined\bysame
\newcommand{\bysame}{\leavevmode\hbox to3em{\hrulefill}\,}
\fi


\begin{thebibliography}{19}

\bibitem{be76}
C.\ Berge, {\em Graphs and Hypergraphs}, Second revised edition, 
North-Holland Mathematical Library 6 (1976) 1-528 (translated from the French by Edward Minieka). 

\bibitem{be89}
C.\ Berge, {\em Hypergraphs.\ Combinatorics of finite sets},
North-Holland Mathematical Library 45 (1989) 1-255 (translated from the French).

\bibitem{bobr97}
M.\ Borowiecki, I.\ Broere, M.\ Frick, P.\ Mih{\'o}k and G.\ Semani{\v s}in,
{\em Survey of hereditary properties of graphs}, 
Discussiones Mathematicae - Graph Theory {\bf 17} (1997) 5--50.\

\bibitem{bomi91}
M.\ Borowiecki and P.\ Mih{\'o}k, {\em He\-re\-di\-ta\-ry pro\-per\-ties of
graphs}, in: V.R.\ Kulli, ed., Advances in Graph Theory (Vishwa International
Publication, Gulbarga, 1991) 42--69.\

\bibitem{brbu99} I.\ Broere and J.\ Bucko, {\em Divisibility in additive
hereditary properties and uniquely partitionable graphs}, Tatra
Mt.\ Math.\ Publ. {\bf 18} (1999), 79--87.\ 

\bibitem{brbu02} I.\ Broere and J.\ Bucko, P. Mih{\'o}k, {\em Criteria for the existence of uniquely partitionable graphs with respect to additive induced-hereditary properties} Discussiones Mathematicae - Graph Theory {\bf 22} (2002), 31--37.\ 

\bibitem{fari02a}
A.\ Farrugia and R.B.\ Richter, Unique factorisation of additive induced-hereditary properties (to appear in Discussiones Mathematicae - Graph Theory).

\bibitem{grhe73}
D.L.\ Greenwell, R.L.\ Hemminger and J.\ Klerlein, {\em Forbidden subgraphs}
Proc.\ 4th S-E Conf.\ Combinatorics, Graph Theory and Computing (Utilitas Math.,
Winnipeg, Man., 1973) 389--394.\

\bibitem{ja02}
J.\ Jakub\'{\i}k, {\em On the lattice of additive hereditary properties of finite graphs}, Discussiones Mathematicae - General Algebra and Applications  {\bf 22} (2002), 73--86.\

\bibitem{mi00} 
P.\ Mih\'{o}k, {\em Unique Factorization Theorem},
Discussiones Mathematicae - Graph Theory  {\bf 20} (2000), 143--153.\

\bibitem{mise00} 
P.\ Mih\'{o}k, G.\ Semani\v{s}in and R.\ Vasky, {\em Additive and hereditary properties of graphs are uniquely factorizable into irreducible factors}, J.\ Graph Theory {\bf 33} (2000), 44--53.\ 

\bibitem{pi91}
B.C.\ Pierce, {\em Basic Category Theory for Computer Scientists}, Foundations of Computing Series, The MIT Press, Cambridge, Massachusetts 1991.\

\bibitem{sc86}
E.R.\ Scheinerman, {\em On the structure of hereditary classes of graphs}, J.\ Graph Theory {\bf 10} (1986) 545--551.\
\end{thebibliography}
\end{document}